\documentclass[12pt]{amsart}
\ifx\pdfoutput\undefined
\usepackage[dvips]{graphicx} 
\else
\usepackage[pdftex]{graphicx} 
\fi
\usepackage{caption}
\usepackage{subfigure}
\usepackage{algorithmic}
\usepackage{algorithm}
\usepackage{amsmath}
\usepackage{amssymb}

\usepackage{color}
\newtheorem{theorem}{Theorem}[section]

\newtheorem{lemma}[theorem]{Remark}

\newcommand{\norm}[1]{\|#1\|}

\newcommand{\coarsevar}{\mathcal{C}}
\newcommand{\finevar}{\mathcal{F}}

\newcommand{\be}{\begin{equation}}
\newcommand{\ee}{\end{equation}}
\newcommand{\beqas}{\begin{eqnarray*}}
\newcommand{\eeqas}{\end{eqnarray*}}
\newcommand{\bea}{\begin{eqnarray}}
\newcommand{\eea}{\end{eqnarray}}
\newcommand{\bc}{\begin{center}}\newcommand{\ec}{\end{center}}
\newcommand{\bi}{\begin{itemize}}
\newcommand{\ei}{\end{itemize}}
\newcommand{\bd}{\begin{description}}
\newcommand{\ed}{\end{description}}
\newcommand{\bdm}{\begin{displaymath}}
\newcommand{\edm}{\end{displaymath}}

\newcommand{\range}[1]{\mbox{range}\left(#1 \right)}

\newcommand{\kernel}[1]{\mbox{null}\left(#1 \right)}

\newcommand{\linspace}[2][ ]{\mathbb{#2}^{#1}}

\newcommand{\mathR}{\mathbb{R}}
\newcommand{\diag}{\mbox{diag}}
\newcommand{\eqnref}[1]{(\ref{#1})}

\newcommand{\AF}[1]{{\color{black}#1}}

\author{M.~Bolten, A.~Brandt, J.~Brannick, A.~Frommer, K.~Kahl, and I.~Livshits}
\begin{document}

\title[Bootstrap AMG for Markov Chains ]{A Bootstrap Algebraic Multilevel method for Markov Chains}
\date{\today}

\maketitle
\begin{abstract}
This work concerns the development of an Algebraic Multilevel
method for computing stationary vectors of Markov chains.
We present an efficient Bootstrap Algebraic Multilevel method for this task.
In our proposed approach, we employ a multilevel eigensolver, with 
interpolation built using ideas based on compatible relaxation, algebraic 
distances, and least squares fitting of test vectors.  Our adaptive variational strategy for computation of the 
state vector of a given Markov chain is then a combination of this multilevel
eigensolver and associated multilevel preconditioned GMRES iterations.
We show that the Bootstrap AMG eigensolver by itself can efficiently compute accurate 
approximations to the state vector.  An additional benefit of the Bootstrap approach is
that it yields an accurate interpolation operator for many other eigenmodes.  This in turn
allows for the use of the resulting AMG hierarchy to accelerate the MLE steps using 
standard multigrid correction steps.  
\textcolor{black}{Further, we mention that our method, unlike other existing multilevel methods for Markov Chains, does not employ
any special processing of the coarse-level systems to ensure that
stochastic properties of the fine-level system are maintained there. }
The proposed approach is applied to a range of test problems, 
involving non-symmetric stochastic M-matrices, showing promising
results for all problems considered.
\end{abstract}

\section{Introduction}\label{sec:intro}
We consider the task of computing a non-zero 
vector $x$ such that
\be
A x = 1\,  x,
\label{eq:ax}
\ee
where $A$ denotes the transition matrix of a given irreducible Markov process
and $x$ is the associated state vector,  an eigenvector of $A$
with eigenvalue equal to one. Since $A$ is irreducible, $x$ is unique up to scalar
factors. The approach considered in this paper is to approximate $x$ 
iteratively using an Algebraic Multigrid
(AMG) method designed to find a non-trivial solution to the eigenvalue problem 
\[
Bx =0 \, x, \enspace \mbox{\AF{where $B = I-A$,}} 
\]
in the setup phase and use the computed AMG method to solve the homogeneous system 
\[ Bx = 0, \].
Our AMG algorithm relies on the Bootstrap 
framework, a fully adaptive algebraic multigrid scheme proposed 
in~\cite{BootAMG2}. For solving 
complex-valued linear systems, this BootAMG (Bootstrap AMG) framework was 
efficiently put into action in~\cite{BootAMG09}. 
In this paper, we develop a variant of BootAMG specifically tailored to 
compute $x$ in (\ref{eq:ax}).
Although we are not (yet) in a position to give a full rigorous 
mathematical analysis of the method, we demonstrate the \AF{efficiency} of the 
BootAMG approach for a series of Markov Chain test problems.

The main new ingredient of our BootAMG method is an adaptive (setup) component 
that simultaneously computes \AF{an} approximation to the state vector and 
Test Vectors (TVs) which are then used in a Least Squares approach to
calculate a multilevel sequence of interpolation operators.
Letting $\ell$ denote a given level of the multilevel hierarchy 
and taking on the finest level $B_0 = B$, 
our adaptive strategy consists of the following steps 
(explained and described in detail in the latter sections):  
\begin{enumerate}
\item relax on the homogeneous system $B_\ell x_\ell  = 0, B_\ell \in \mathR^{n_\ell \times n_\ell}$ to improve
the set of test vectors,  $x^{(1)}_{\ell}, ..., x^{(r)}_{\ell}$, approximating the state vector; $n_\ell$ here is the problem size on level $\ell$;\item compute a sparse interpolation operator $P_\ell \in \mathR^{n_{\ell+1} \times n_\ell}$ using a least-squares fit based on $x^{(1)}_{\ell}, ..., x^{(r)}_{\ell}$; 
\item compute the system matrix $B_{\ell+1} \in \mathR^{n_{\ell+1} \times n_{\ell+1}}$ on the next coarser level 
      as $B_{\ell+1} = Q_\ell B_\ell P_\ell$  via a Petrov-Galerkin
      approach, with the restriction $Q_l \in \mathR^{n_{\ell} \times n_{\ell+1}}$ 
      representing a simple averaging operator; also calculate the mass matrix $T_{l+1}  = Q_\ell T_\ell P_\ell$ (with $T_0 = I$) for the next leg of BootAMG;  
\item obtain the coarse-level series of test vectors 
  $$ x^{(j)}_{\ell+1}= R_\ell x^{(j)}_{l},\ j = 1,...,r,\quad \text{with\ } 
    R_\ell \in \mathbb{R}^{n_{l+1} \times n_l}\quad \text{being injection}.$$ 
\end{enumerate}
This process is applied recursively until a level with sufficiently small problem size 
(the coarsest level, $\ell = L$) is reached.  
There,  an exact eigensolve for 
\[ B_\ell x = \lambda T_\ell x \]
is performed and then used to improve the given approximations to the state vector 
on increasingly finer levels.  Along with the state vector's approximations $x^{(0)}_\ell$,  
some (small number) of near kernel eigenfunctions $(\lambda_i, x^{(i)}_\ell)$, from the coarsest level
are interpolated
to increasingly finer levels, where they are  relaxed and then used   as initial guesses to the corresponding  fine-level eigenvectors. On each such level, the relaxation is performed on the \AF{homogeneous} system 
\[ (B_\ell - \lambda_i T_\ell) x^{(i)}_\ell = 0, \]
and the approximations to $\lambda_i$ are updated  (except for $\lambda_0$ which is known to be exactly zero.)  Each $x^{(i)}_\ell$ is already a reasonable approximation to \AF{an} eigenfunction on level $\ell$, i.e, 
$$ B_\ell x^{(i)}_{\ell} \approx \lambda_i^\ell  T_\ell x^{(i)}_\ell $$
for \AF{an} $\ell^{th}$ level eigenvalue $\lambda_i^\ell$. Therefore, a good eigenvalue approximation can be obtained by 
\[ \hat{\lambda}_i = \displaystyle\frac{(B_\ell x^{(i)}_{\ell}, v)}
{(T_\ell x^{(i)}_{\ell}, v)}
\] for any vector $v$. We choose $v = x^{(i)}_\ell$ and thus calculate an approximation to eigenvalue using what would be a Rayleigh quotient for Hermitian matrices. The change in eigenvalue  approximation  for  each  vector $x^{(i)}, \,  i > 1$
 after 
relaxation provides a measure of their accuracy -- if after relaxation the relative change is significant, 
we add  $x^{(i)}$ to  the set of test vectors.  
Otherwise, these eigenvectors are well approximated by the existing  interpolation operator to the current level and thus
need not be built into the related coarser space. The state  eigenvector approximation is always included in the test set.  After all test sets on all levels, including the finest,  are updated,  all  $P_\ell$,  $B_\ell$, and $T_\ell$  are recalculated, and the setup stage is complete.
 
Once the adaptive setup is finished, we use the resulting V-cycle AMG method as
a preconditioner to (full) GMRES, applied to the homogeneous system $Bx = 0$. 
This iteration converges rapidly to the state vector, since as we demonstrate later, the BootAMG preconditioner
efficiently separates the state vector from the field of values  of $B$ on an
appropriate complementary subspace.  This use of both MLE steps and 
V-cycle preconditioned GMRES iterations is the other main new idea in the present paper.

The remainder of this paper is organized as follows.
Section~\ref{Sec:MCO} contains an introduction to Markov chain systems 
and a general review of AMG approaches 
for computing state vectors of Markov chains.  
Section~\ref{sec:bootamg} contains a general
description of the Bootstrap AMG ideas along with full algorithmic descriptions
of realizations that we use in this paper. In Section~\ref{sec:Alg} we 
discuss our BootAMG algorithm for Markov chains and theory for our AMG-GMRES 
solver.  Numerical experiments are presented in Section \ref{sec:NumRes} and concluding remarks are given in Section \ref{sec:Conclusions}.

\section{Overview of Solvers for Markov Chain Systems}\label{Sec:MCO}
The transition matrix  
$A \in \mathR^{n \times n}$ of a Markov process contains as its entries $a_{ij}$, 
the transition probabilites from state $i$ to state $j$,  $a_{ij} \geq 0$
for all $i,j$.  Matrix $A$ is column stochastic, i.e.\ $A^t\mathbf{1}
= \mathbf{1}$, with $\mathbf{1}$ being the vector
of all ones. It is always possible to eliminate self-transitions~\cite{SE81}, so we assume 
$a_{ii} = 0$ for all $i$ from now on. The state vector   
$x$ satisfies 
$$A x = x$$
with $x \neq 0$, $ 0 \leq x_i $. By the Perron-Frobenius theorem 
(cf. \cite[p.~27]{Berman_Plemmons_1994}, e.g.) such a vector $x$ always exists.

A general square matrix $A \in \mathR^{n \times n}$ induces a directed graph 
$G(A) = (V,E)$ with vertices $V=\{1,\ldots,n\}$ and directed edges 
$E=\{ (i,j) \subset V^2: i \neq j \mbox{ and } a_{ij} \neq 0\}$. 
Two vertices $i$ and $j$ in $G(A)$ are said to be strongly
  connected (as in graph theory)
if there exist directed paths in $G(A)$ from $i$ to $j$ and from $j$ to $i$.
Since this is an equivalence relation on the vertices, it induces a partitioning 
of $V$ into the strong components of $G(A)$. If $G(A)$ has exactly one strong
component, the matrix $A$ is called irreducible; otherwise it is called reducible
which is precisely the case when there exists a permutation matrix $\mathcal{P}$ such that 
$$\mathcal{P}^tA\mathcal{P} = 
\begin{pmatrix}
  A_{11} & 0 \\
  A_{12} & A_{22} 
\end{pmatrix}, 
$$
where  $A_{11}$ and $ A_{22}$  are  square submatrices.

If $A$ is irreducible --- which we assume throughout --- the Perron-Frobenius 
Theorem guarantees that the state vector $x$ has all positive values and
is unique up to a scalar factor.  
From this theorem it also follows that the spectral radius of $A$ satisfies 
$\rho(A) = 1$, implying that $B = I -A$ is a singular $M$-matrix; recall that $ 0 \leq a_{ij} \leq 1$.  
In addition, $B$ is irreducible, because of the irreducibility of $A$.

The idea of using Multigrid (MG) to compute the state vector of an irreducible transition matrix 
is not new.  Numerous approaches have been explored in 
the past, in both the aggregation-based~\cite{HDeSterck_08,HDeSterck_09,HSIMON_61,YTakahashi_75} 
and Classical~\cite{Evirnik_07} MG settings.  In~\cite{HSIMON_61}, 
the idea of using a partition of 
unity type aggregation method was first explored.  Later, in~\cite{HDeSterck_08}, 
an extension of this approach based on Smoothed Aggregation
was developed for application to the page rank problem and further generalized 
in~\cite{HDeSterck_09}.  For an overview of Classical MG type methods 
(and preconditioners) we direct the reader to~\cite{Evirnik_07} and the references therein.  

Our MG solver can be loosely categorized as a variation of the Classical AMG 
approach in the following way.  We choose  the coarse level
variables as a subset of the fine level variables.  In this aspect, 
our approach is most closely related to the recent work by Virnik~\cite{Evirnik_07};   we too employ the AMG-type method as a preconditioner for GMRES.   
Our proposed new approach, however, differs from all previous AMG-type solvers
for Markov chains in several ways. Most importantly, we build interpolation 
adaptively using a Least Squares approach, and  we combine the multilevel eigensolver with application of correction steps using an AMG-GMRES preconditioner. Unlike other methods, we thus (experimentally) observe an excellent convergence without the need to  recompute 
interpolation and the entire  coarse-level hierarchy at each iteration.
In addition, our  solver yields an efficient and straightforward 
method for Markov chains which computes the state vector to 
arbitrary accuracy without the need for {\em lumping} (see~\cite{HDeSterck_09}) and other such approaches 
typically employed to maintain the
 stochastic properties of the transition matrix and state vectors on all coarse 
levels of the MG hierarchy.

\section{Bootstrap AMG}\label{sec:bootamg}
In this section, we provide a general review of the
Bootstrap AMG framework~\cite{BootAMG2,ABIL09} for building 
MG interpolation together with 
some heuristic motivation of  our choices for the 
individual components of the BootAMG multigrid algorithm.

The first key ingredient to any AMG method is its {\em relaxation} or 
{\em smoothing} iteration. For a given system matrix $B$ it is a splitting
based iteration 
\begin{equation} \label{amg_relax:eq}
   Mu^{\nu+1} = Nu^{\nu} + b, \enspace \mbox{ where } B = M-N, \, M \mbox{ non-singular}, \,\nu=0,1,\ldots
\end{equation}
with its ``error propagation matrix'' given by $M^{-1}N = I -M^{-1}B$. 
Here, $M$ is
chosen such that after a few iterations the error $e^{\nu}$ w.r.t.\ the solution 
of $Bx = b$ is {\em algebraically smooth},
i.e.\ $\|Be^{\nu} \| \ll \|e^{\nu}\|$. In many situations this can be achieved by point-wise
Gauss-Seidel or (under-relaxed) Jacobi iterations. Indeed, in our Markov chain
setting we used $\omega$-Jacobi relaxation with $\omega = \AF{0.7}$ on all levels, 
i.e., we took $M = \frac{1}{\omega} \diag(B)$, resulting in
$I - M^{-1}B = I - \omega \cdot (\diag(B))^{-1}B$. 

(In case a diagonal entry $b_{ii}$  of $B$ happens to be zero for some row(s)
$i$,  Kaczmarz or some other distributive relaxation scheme can be applied to the $i^{th}$ equation.) \\

The Bootstrap AMG process can now be described as follows.
Coarse variables are selected as a subset of the fine variables using a full coarsening (for problems on structured grids) or, otherwise,  compatible relaxation (CR) coarsening 
scheme \cite{ABrandt_2000}. The CR scheme can be either started from scratch, or, 
if geometric information is given and a suitable candidate set of coarse variables is
known,  such set can be tested and improved by CR.
Interpolation is then computed using a Least Squares based approach.
We mention  that once a tentative Multigrid hierarchy
has been defined, it can be used to further enhance the set of test vectors.

\subsection{Choosing the coarse variables: compatible relaxation}\label{sec:CompRel}
In the AMG setting, CR is a relaxation-based 
coarsening process which can be viewed as a special case
of a general approach for constructing coarse-level descriptions
of given fine-level systems, including non-stationary, highly
non-linear, and also non-deterministic systems~\cite{ABrandt_SysUpsc08}.
The basic idea of CR is to use the given relaxation scheme \eqnref{amg_relax:eq}, 
restricted to appropriately defined subspaces, to measure the quality
of the given coarse space and also to iteratively improve it if needed.
We proceed with a brief overview of CR and its use in AMG coarsening.
A detailed discussion, theory and comparisons between various
measures of the quality of coarse spaces and their relations to
compatible relaxation schemes are presented
in~\cite{ABrandt_2000, JBrannick_Thesis, JBrannick_RFalgout, 
PanayotRob_2003, FVZ, Olivne}.

\subsubsection{Classical AMG CR-based coarsening}
Assume that the set of coarse-level variables, $\mathcal{C}$, is a
subset of the set of fine-level variables, $\Omega$.  Under this assumption, 
one possible form of CR is given by $\mathcal{F}$-relaxation for the
homogeneous system --- relaxation applied only to the set 
of $\mathcal{F}$ variables, with $\mathcal{F} := \Omega \setminus \mathcal{C}$.
Given the partitioning of $\Omega$ into $\mathcal{F}$ and $\mathcal{C}$,
we have
$$
u = \begin{pmatrix} u_f \\ u_c \end{pmatrix}, \enspace 
B = 
\begin{pmatrix}
    B_{ff}  & B_{fc}   \\
    B_{cf} & B_{cc} 
\end{pmatrix},
\enspace \mbox{and} \enspace
M = 
\begin{pmatrix}
    M_{ff}  & M_{fc}   \\
    M_{cf} & M_{cc} 
\end{pmatrix},
$$
assuming the equations are permuted such that the unknowns in $\mathcal{F}$ 
come before those in $\mathcal{C}$. 
The $\mathcal{F}$-relaxation of CR is then defined by 
\begin{equation} \label{frelax:eq}
   u_f^{\nu+1} = (I-M_{ff}^{-1}B_{ff})u_f^{\nu} = E_{f}u_{f}^{\nu}. 
\end{equation}
If $M$ is symmetric, the asymptotic convergence rate of CR 
\[
 \rho_{f}
=\rho(E_{f}), 
\]
where $\rho$ denotes the spectral radius, provides a measure of the quality of the
coarse space, that is, a measure of the ability of the set of coarse
variables to represent error not eliminated by the given fine-level relaxation.
This measure can be approximated using $\mathcal{F}$-relaxation for the homogeneous 
system with a random initial guess $u_f^{0}$. Since $\lim_{\nu \to \infty} 
\|E_f^\nu\|^{1/\nu} = \rho(E_f)$ for any norm $\| \cdot \|$, the measure
\begin{equation} \label{meas1:eq}
   \left(\|u_f^\nu\| / \|u_f^0\| \right)^{1/\nu}
\end{equation}
estimates $\rho_f$. 
  
In choosing $\mathcal{C}$, we use the CR-based coarsening 
algorithm developed in~\cite{JBrannick_RFalgout}.  This 
approach is described in Algorithm \ref{alg:CR}.

\begin{algorithm}
\caption{compatible\_relaxation \hfill \{Computes $\mathcal{C}$ using Compatible Relaxation\}\label{alg:CR}}
\begin{algorithmic}
\STATE \textit{Input:} $\mathcal{C}_{0}$ \COMMENT{$\mathcal{C}_{0} = \emptyset$ is allowed}
\STATE \textit{Output:} $\mathcal{C}$
\STATE Initialize $\mathcal{C} = \mathcal{C}_{0}$
\STATE Initialize $\mathcal{N} = \Omega \setminus \mathcal{C}$
\STATE Perform $\nu$ CR iterations \eqnref{frelax:eq} with components of $u^0$
       randomly generated
\WHILE{$\rho_{f} > \theta$}
\STATE $\mathcal{N} = \{i \in \Omega \setminus \mathcal{C}: \frac{\displaystyle 
|u_{i}^{\nu}|}{ \displaystyle |u_{i}^{\nu-1}|} > \theta \}$
\STATE $\mathcal{C} = \mathcal{C} \cup \{\text{independent set of\ }  \mathcal{N}\}$
\STATE Perform $\nu$ CR iterations \eqnref{frelax:eq} with components of $u^0$
       randomly generated
\ENDWHILE
\end{algorithmic}
\end{algorithm}

In our numerical experiments for Markov chains, we use weighted Jacobi CR, set the CR tolerance $\theta = .85$, the number of CR sweeps $\nu =8$, choose the components of $u^0$ to be uniformly distributed in the interval $[1,2]$, and select $\mathcal{C}_{0}$
using the standard independent set algorithm (see
\cite{oosterlee_book}) based on the full graph of the system
matrix. For further information on CR we refer the reader to~\cite{JBrannick_RFalgout}.

\subsection{Building Bootstrap AMG interpolation}\label{sec:BootAMG}
We now outline the Least Squares (LS) approach for defining
interpolation, see also \cite{BootAMG09}. We assume that 
the set ${\mathcal C}$ of coarse variables is given, i.e., it has been  
previously determined, for instance,  by geometric coarsening and/or using compatible 
relaxation.  In the LS approach, the interpolation operator $P$ is chosen to approximate a given (specifically chosen) 
set of test vectors. We define $c$ as the maximum 
number of coarse-level variables used to interpolate to a single
fine-level variable, or equivalently the maximum number of non-zero
entries in any row of $P$. 
The key ingredient of the BootAMG setup
lies in the use of several test vectors (TVs) that collectively 
should represent 
those error components not reduced by relaxation.
We assume for now that such a set of test vectors, $\mathcal{U} :=\{x^{(k)}\}_{k = 0}^r$ 
is known on some fine level. 
The rows of the prolongation operator $P$ are then obtained individually.
For each variable $i \in \mathcal{F}$, we first determine a set of its neighboring coarse-level  variables, 
$\mathcal{N}_i$,  using  the directed graph $G(A)$ 
\begin{equation} \label{eq:neighbor_def}
   \mathcal{N}_i = \{ j \in \mathcal{C}: \mbox{ there is a path of length $\leq z$ in $G(A)$ from $i$ to $j$ } \},
\end{equation} where $z \leq 3$, such that $\mathcal{N}_{i}$ is a
subset of the local graph neighborhood of $i$.
We then determine an appropriate set of (coarse level) interpolatory variables
 $\mathcal{J}_{i}
\subseteq \mathcal{N}_i$  with $| \mathcal{J}_i | \leq c$ that yields the best fit of the LS interpolation for point $i \in \mathcal{F}$. \AF{For ease of presentation let us
use the indices $j \in \mathcal{J}_i$ of the fine level variables to address the columns of $P$. We then } define the local least squares functional  for  the nonzero entries
$P_i = \{p_{ij}: j \in \mathcal{J}_i\}$ of row $i$ of $P$ as
\begin{equation}\label{eq:LSfuncrowi}
L (P_i;\mathcal{J}_i) = 
\sum_{k=0}^r\omega_k \left( x_{i}^{(k)} - \sum_{j\in \mathcal{J}_i} p_{ij} \AF{x}_{j}^{(k)} \right)^{2} . 
\end{equation}

The task is then to find a set $\mathcal{J}_i$ of interpolating 
points for which 
the minimum of $L$ is small and \AF{to obtain} the corresponding values $p_{ij}$ of 
the minimizer that  yield the coefficients for the interpolation operator. 
Generally, the weights $\omega_{k}$ should be chosen according to the 
algebraic smoothness of $x^{(k)}$  to bias the least squares functional 
towards the smoothest vectors. We do so by using the square of the
inverse of the \AF{norm of the} residual of the \AF{homogeneous} problem $B_\ell x^{(k)} = 0$
for each $x^{(k)}$.
Here $\ell$ is the current level to which operator $P$ interpolates
from the next coarser level $\ell+1$. 

\AF{Our approach for solving this task is given as Algorithm~\ref{alg:LSIS}. It is based on
a {\em greedy strategy} to find an appropriate set of interpolatory points $\mathcal{J}_i$ for 
each fine point $i$; see \cite{BootAMG09} and \cite{kahl_thesis:2009}.}


\bigskip
\begin{algorithm}
\caption{ls\_interpolation \hfill \{Computes Least Squares based interpolation\} \label{alg:LSIS}}
\begin{algorithmic}
  \STATE \textit{Input:} $\mathcal{U}$, $c$
  \STATE \textit{Output:} interpolation $P$
  \STATE $X^{(k)} = R \,x^{(k)}, x^{(k)} \in \mathcal{U}$
\FOR{$i \in \mathcal{F}$}
	\STATE Take $\mathcal{N}_{i}$ from \eqnref{eq:neighbor_def}
	\STATE Set $\mathcal{J}_{i} = \emptyset$
        \REPEAT
        \STATE determine $g^{*} \in \mathcal{N}_i$ s.t.\
         $\underset{P_i}{\min} \, L(P_i ; \mathcal{J}_i \cup \{g^*\}) = \underset{g \in \mathcal{N}_{i}}{\min}\, \underset{P_i}{\min} \, L(P_i,\mathcal{J}_i \cup \{g\})$		
        \STATE Set $\mathcal{N}_{i} = \mathcal{N}_{i} \setminus \{g^{*}\}$ and $\mathcal{J}_{i} = \mathcal{J}_{i} \cup \{g^{*}\}$
	\UNTIL{$|\mathcal{J}_{i}| \geq c$ or $\mathcal{N}_{i} = \emptyset$}
\ENDFOR
\end{algorithmic}
\end{algorithm}

\section{MLE and AMG-GMRES Algorithms}\label{sec:Alg}
In this section, we discuss the main ideas and their implementations
in our BootAMG-based algorithm for 
computing non-trivial solutions $x$ to $Bx = 0$, where as before $B = I - A$,  $A$ 
is an irreducible column stochastic transition matrix so that $B$ is a non-symmetric 
and singular $M$-matrix. 
The goal of each MLE step is twofold: (1) to compute approximations to
the solution $x$ and also to compute suitable test vectors that \AF{allow}
us to improve the accuracy of Least Squares interpolation
and (2) to update the multigrid hierarchy, used to precondition GMRES. 

Such a hierarchy consists of a sequence of prolongation and restriction operators 
and systems of coarse-level equations. The coarse-level equations are defined 
using the Petrov-Galerkin approach. \AF{More precisely,} $B_{l+1} =
Q_lB_lP_l$, with restriction $Q_{l}$ given by the averaging operator. 
\AF{This means that each column of $Q_l$ has identical nonzero entries at exactly the positions corresponding to $\mathcal{J}_i$ and that $\mathbf{1}^tQ_{l} = \mathbf{1}^t$.  
This implies $B_{l}^{t}\mathbf{1} = \mathbf{0}$ for all levels $l$.} 

\subsection{Multilevel Eigensolver (MLE)}\label{sec:RecSCS}
Our approach for building a multilevel method relies on the BootAMG framework. 
We construct a sequence of increasingly coarser-level descriptions of the given fine-level
Markov Chain system using the LS-based fit discussed above.
The coarse-level test vectors are calculated from the fine-level test set as 
$X^{(k)} = Rx^{(k)}$. 
On the \AF{coarsest level}, the eigenproblem is solved exactly 
for some small number of the lowest (in absolute value)  eigenvectors, which are then directly  interpolated and relaxed on finer levels, becoming increasingly more accurate approximations to the finest-level eigenfunctions (such treatment is similar to the Exact Interpolation Scheme, e.g, \cite{IL08}).  These 
eigenfunction approximations are added to test sets on each level.  
Note that the  initial (random) test functions are not relaxed during this stage -- otherwise they may become almost linearly dependent which may result in \AF{ill-posed  local} least squares problem.
With the  improved  TV sets (enriched by eigenfunction approximations),    
new prolongation operators are then calculated.

We note that the setup cycle in Algorithm~\ref{alg:BootAMGWSC}
also yields an approximation to the state vector which we use as an initial guess in our AMG-GMRES iterations.  This approach can also be used as a stand-alone solver for these systems although, in general, it will be far more expensive than using several AMG-GMRES steps
in addition to the MLE.  We note that for certain problems (e.g., when several eigenvalues of $B$ are close to zero) and for larger problem sizes, 
it may be necessary to alternate several times between the MLE and AMG-GMRES iterations to find an accurate approximation of $x$.  
In our tests,  we found it sufficient to use either a traditional $V$-cycle ($\mu =1$) or $W$-cycle ($\mu = 2$)
for the MLE step followed by a small number of AMG-GMRES steps to compute an approximation to the state vector.
Ultimately, in our MLE approach, the sequence of prolongation operators 
needs to be accurate for only the smoothest error component -- the kernel 
of the finest-level operator $B = I-A$. 
The accurate resolution of many other near-kernel components by our LS-based interpolation operator
allows for \AF{fast convergence of} 
MG preconditioned GMRES,
whereas the kernel component of $B$ is the state vector we aim to compute. 
Algorithm \ref{alg:BootAMGWSC} contains a pseudo-code of our implementation of the MLE iteration. 

\begin{algorithm}
\caption{bootamg\_mle \hfill \{Implements the Bootstrap AMG MLE scheme\}\label{alg:BootAMGWSC}}
\begin{algorithmic}
  \STATE \textit{Input:} $B_{l}$ ($B_0 = B$), $T_{l}$, ($T_{0} = I$),  $\mathcal{U}_{l}$ (set of test vectors on level $l$)
  \STATE \textit{Output:} updated $\mathcal{U}_{l}$ and $(\Lambda_{l}, \mathcal{V}_{l})$, approximations to the lowest eigenpairs; 
  \IF{$l = L$}
    \STATE Compute $\mathcal{V}_{L} = \{x^{(i)}\}_{1,\ldots,k}$,
    s.t.\ $B_{L}x^{(i)} = \lambda^{(i)} T_{L}x^{(i)},
    |\lambda^{(1)}| \leq |\lambda^{(2)}| \leq \ldots \leq |\lambda^{(k)}|$.
  \ELSE
    \STATE Relax $B_{l}x^{(j)} = 0,\ x^{(j)} \in \mathcal{U}_l$
    \STATE Set $\mathcal{V}_{l} = \emptyset$
    \FOR{$m = 1,\ldots,\mu$}
      \STATE $P_{l} = \mbox{ls\_interpolation}(\mathcal{U}_{l}\cup \mathcal{V}_{l}, c)$
      \STATE Compute averaging $Q_{l}$ with $\mbox{sparsity}(Q_{l}) = \mbox{sparsity}(P_{l})$
      \STATE $B_{l+1} = Q_{l} B_l P_{l}$
      \STATE $T_{l+1} = Q_{l} T_l P_{l}$
      \STATE $\mathcal{U}_{l+1} = \{ R_{l}x, x \in \mathcal{U}_l\}$ 
      \STATE $\mathcal{V}_{l+1} = \mbox{bootamg\_mle}(B_{l+1},T_{l+1},\mathcal{U}_{l+1})$
      \STATE ${\mathcal{V}}_l = \{ P_{l}x, x \in \mathcal{V}_{l+1}\}$
      \FOR{$i = 1,\ldots,|{\mathcal{V}}_l|$}
         \STATE Relax $(B_{l} - \lambda_i T_{l} )x^{(i)} = 0$, $x^{(i)} \in \mathcal{V}_l$
	 \STATE Update ${\lambda}_i = \displaystyle
	 \frac{\langle B_{l} x^{(i)}, x^{(i)}\rangle}{\langle T_{l} x^{(i)}, x^{(i)}\rangle}$
      \ENDFOR
    \ENDFOR
   \ENDIF
\end{algorithmic}
\end{algorithm}

\begin{figure}[!ht]
  \centering 
  \scalebox{.65}{\input{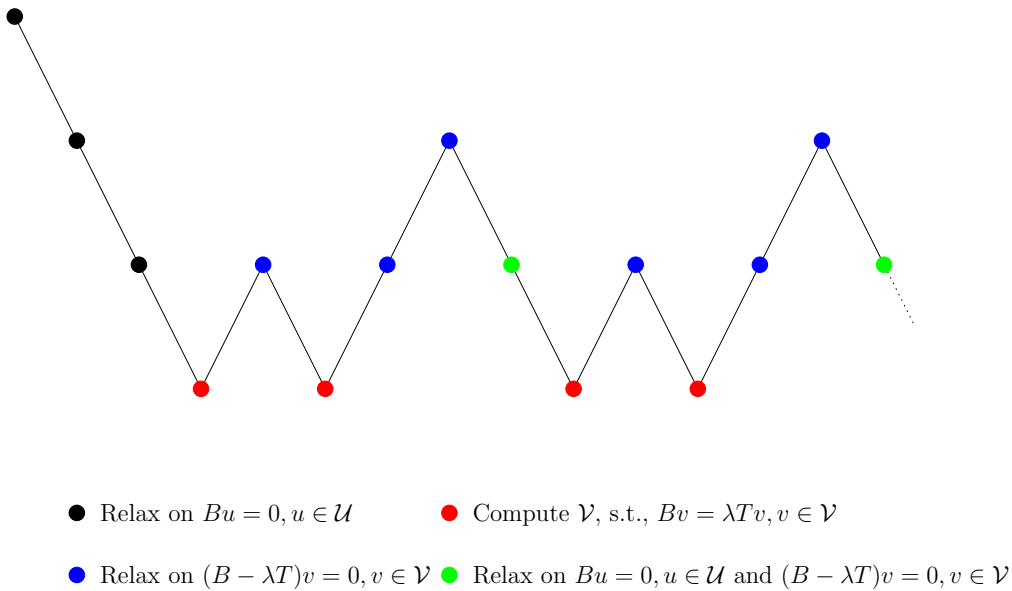}}
  \caption{Bootstrap AMG setup W-cycle.\label{fig:boot:setupcycle}} 
\end{figure}

In Figure~\ref{fig:boot:setupcycle} a possible setup cycle is visualized.  We note that
at \AF{each square} one has to decide whether to recompute $P$ or advance to the next
finer grid in the multigrid hierarchy.
The illustrated cycle resembles a $W$-cycle, but any other cycling strategy can 
be applied to the setup. 

\subsection{Multigrid preconditioned GMRES} 
The multilevel eigensolve steps described above yield increasingly 
better approximations to the state vector and other vectors
that cannot be removed efficiently by the Multigrid relaxation.

Although we are only interested in computing the state vector, the LS-based 
Multigrid hierarchy is able to resolve a larger subspace. This
leads to the idea of exploiting this richness of the given hierarchy
for use in MG correction steps -- in addition to the discussed MLE steps.

To illustrate the effect of MG correction steps applied to
the homogeneous problem $Bx = 0$, we start by analyzing simple relaxation
schemes for the steady state problem, $Ax = x$.
As the steady state solution is the eigenvector corresponding to the eigenvalue
with largest absolute magnitude, a power iteration
\begin{equation}\label{eq:powerit}
  x^{k+1} = Ax^{k}
\end{equation} is guaranteed to converge to the solution. However, convergence can be slow
if $A$ has other eigenvalues with absolut value close to one.

Such power iterations applied to the steady state problem are in turn equivalent 
to applying a Richardson iteration to the homogeneous system $\left(I - A\right)x = 0$.
A natural modification, facilitating that the field of values of $A$ is contained in
the unit circle, is then given by a suitable under-relaxed iteration, yielding the
error propagator
\begin{equation}\label{eq:omegarich}
  e^{k+1} = \left(I - \tau B\right)e^{k}.,
\end{equation} 
\AF{which translates into an identical relation for the iterates, 
\begin{equation}\label{eq:omegapower}
  x^{k+1} = \left(I - \tau (I-A) \right)x^{k},
\end{equation}}
 so that the iteration can be interpreted as a modified power method.

\begin{figure}
  \subfigure[Spectrum and field of values for Richardson error propagator.\label{fig:specplotJ}]{\includegraphics[width=.3\textwidth]{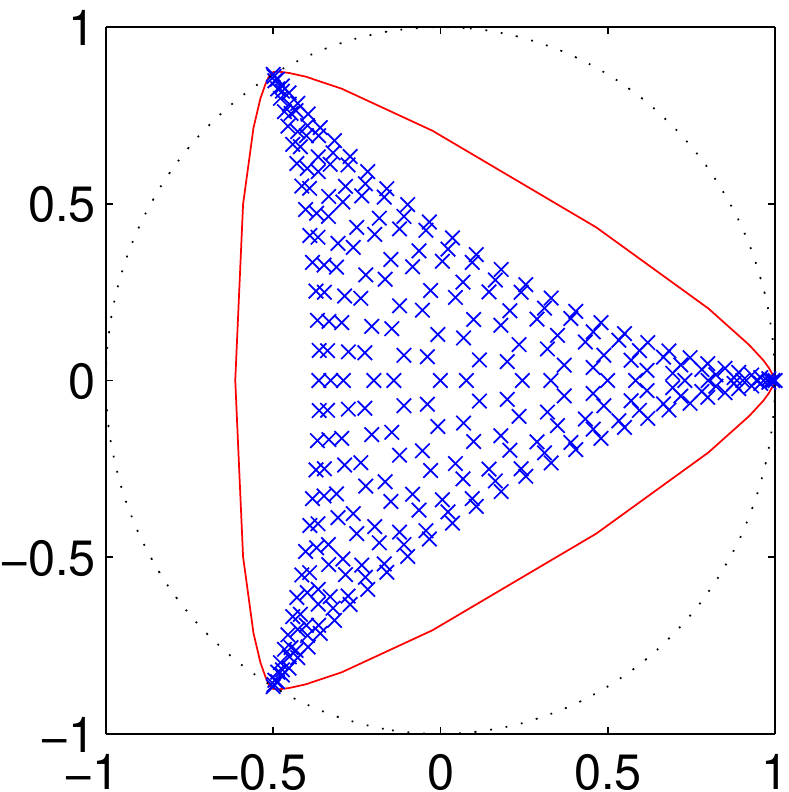}}\hfill
  \subfigure[Spectrum and field of values for $\tau$-Richardson.\label{fig:specplotwJ}]{\includegraphics[width=.3\textwidth]{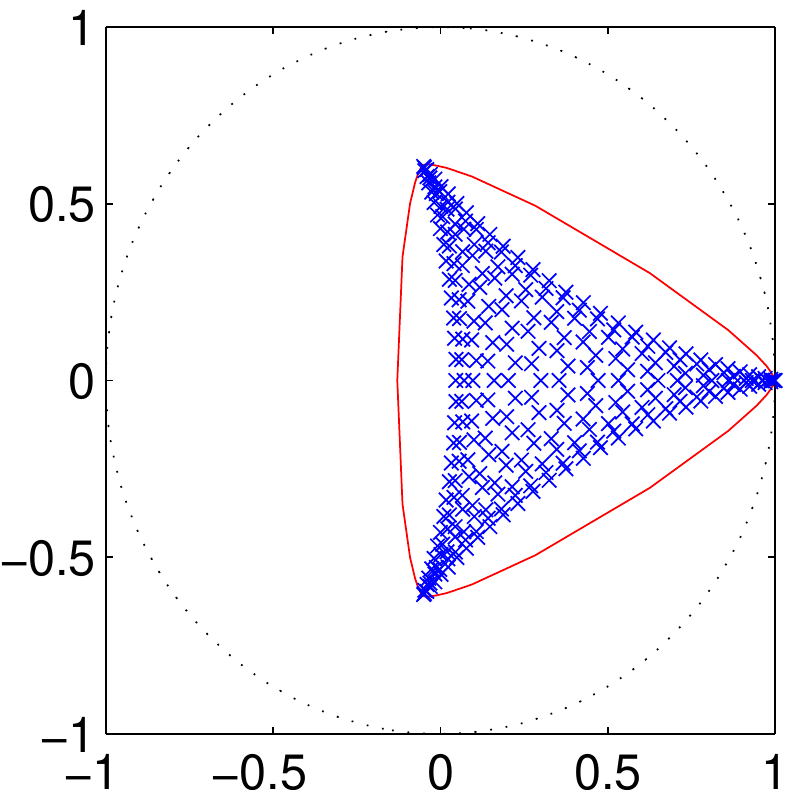}}\hfill
  \subfigure[Spectrum and field of values of Multigrid error propagator.\label{fig:specplotmgJ}]{\includegraphics[width=.3\textwidth]{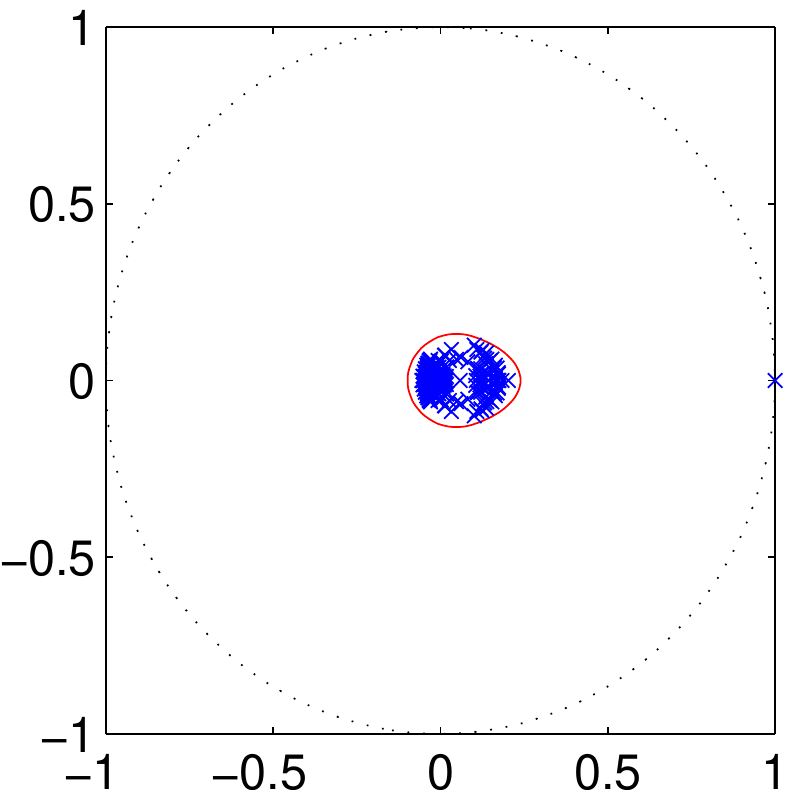}}
  \caption{Spectra and field of values of the Richardson, $\tau$-Richardson and Multigrid error propagators for two-dimensional tandom queuing problem with $n = 33^2$.\label{fig:specplot}}
\end{figure}

In Figures~\ref{fig:specplotJ} and~\ref{fig:specplotwJ} the spectra of $A$ and 
of $I - \tau\left(I - A\right)$ are depicted for a characteristic two-dimensional test problem, 
along with the field of values.
\AF{The error propagator of the MG V-cycle applied to
$Bx = 0$ is given by
\begin{equation}\label{eq:mgerrorprop}
  E = \left(I - MB\right)\left(I - P{E_{c}}QB\right)\left(I - MB\right), 
\end{equation} where ${E_{c}}$ denotes the error propagation operator on the
next coarser level. This recursive definition terminates at the coarsest level with
$E = B^\dagger$, the Moore-Penrose inverse of the matrix on the coarsest level, since we 
compute the minimal norm solution of the respective system at that level.}
The error
propagator can be rewritten as
\begin{equation}\label{eq:mgpower}
  \AF{I - \mathcal{C}B},
\end{equation} 
\AF{so that comparing with the one from \eqnref{eq:omegarich} 
(relaxed Richardson) we see that we can interprete the multigrid V-cycle as yet another 
modification of the power iteration.}
In Figure~\ref{fig:specplotmgJ}, the spectrum and field of values of the MG
V($2,2$)-cycle preconditioned matrix \AF{(\ref{eq:mgpower})} is depicted.  For this case, it is clear that applying
power iterations to this preconditioned operator will converge rapidly to the steady state vector.

To further accelerate this approach, instead of the straightforward power method,  we consider  MG preconditioned GMRES steps.  As we show below, such
iteration is guaranteed to converge to the solution of our Markov
Chain systems. 

\begin{lemma}\label{markov:rangenullrelation}
Let $A$ be a column-stochastic irreducable operator and $B = I - A$.
Then we have
\begin{equation}\label{eq:markov:rangenullrelation}
  \range{B} \cap \kernel{B} = \{0\}.
\end{equation}
\begin{proof}
  As $A$ is column-stochastic we know that $A^{t}\mathbf{1} = \mathbf{1}$.
  Hence $B^{t}\mathbf{1} = 0$. Furthermore we know from the Perron-Frobenius
  theorem that $Bx = 0$ is uniquely solvable up to a scalar factor and
  the solution $x^{*}$ is strictly positive (component-wise).
  Thus, as for all $y \in \range{B}$ there exists a $z$ such that \AF{$y = Bz$} 
  we have $\langle \mathbf{1}, y \rangle = \langle \mathbf{1}, \AF{Bz} \rangle = \langle \AF{B^{t}}\mathbf{1}, z \rangle = 0$.
  With this and $\langle \mathbf{1}, x^{*} \rangle \neq 0$ as $x^{*}$ is 
strictly positive we get~\eqref{eq:markov:rangenullrelation}.
\end{proof}
\end{lemma}

In~\cite[Theorem 2.8]{KHayami_MSugihara_2009} it is shown that
GMRES determines a solution of $Bx = b$ for all 
$b \in \range{B}, x_{0} \in \linspace[m]{R}$ iff 
$\range{B} \cap \kernel{B} = \{0\}$. Due to 
Lemma~\ref{markov:rangenullrelation} the assumptions of this theorem
are fulfilled for $B = I - A$, where $A$ is a column-stochastic 
irreducible \AF{operator} arising in Markov-Chain models. \AF{ 
There seems to be no simple way to extend Lemma~\ref{markov:rangenullrelation} to
the multigrid preconditioned matrix. So the observation that 
GMRES for the multigrid preconditioned matrix worked out very well in all our experiments is
not backed by theory so far.}
%
\subsection{A multigrid preconditioned Arnoldi method}
\AF{Instead of solving the homogeneous linear system with preconditioned GMRES, we can also 
perform the Arnoldi method for the preconditioned matrix $\mathcal{C}B$. We thus compute 
orthonormal vectors $v_1,v_2,\ldots$ in the usual way such that $v_1,\ldots,v_k$ is a basis
of the Krylov subspace $\mbox{span}\{y,\mathcal{C}y,\ldots,(\mathcal{C}B)^{k-1}y\}$. Denoting
$H_k \in \mathR^{k \times k}$ the orthogonal projection of $\mathcal{C}B$ onto the Krylov subspace,
i.e. $H_k = V_k^t ( \mathcal{C}B) V_k$ with $V_k = [v_1 | \ldots | v_k]$, we compute the eigenvalue
closest to 0 of $H_k$ and the corresponding eigenvector $\eta$. We then take the Ritz vector $V_k \eta$
as an approximation for the steady state vector. For the vector $y$ upon which we build the Kryolv subspace we take
the approximation $x$ for the steady state vector resulting from the MLE setup phase.}   

\section{Numerical Results}\label{sec:NumRes}
In this section, we provide results obtained using our Bootstrap AMG 
method when applied to a series of Markov Chains.
Our numerical tests consist of our BootAMG-based approach to three Markov chain models that
can be represented by planar graphs. Each of the models has certain
interesting characteristics that pose problems for the solver. 

\subsection{Test Problems}
\begin{figure}[!b]
  \centering 
  \scalebox{.75}{\input{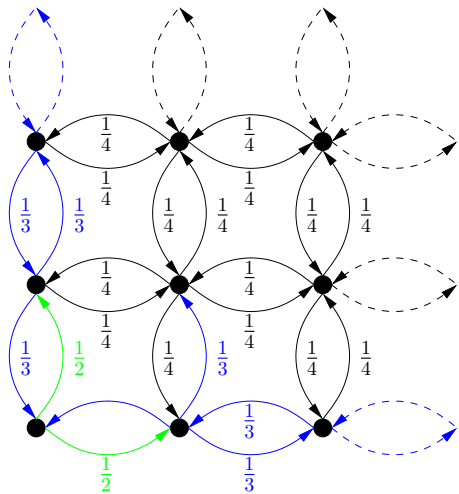}}
  \caption[Uniform network model]{Uniform network model on a two-dimensional equidistant grid.\label{fig:nr:scalar:uniformnet}}
\end{figure}

\AF{We begin our experiments with a very simple model.}
The uniform two-dimensional network can be seen as the Markov chain
analogue of the Laplace operator. It is defined on an equidistant grid
$\Omega$ of size $N \times N$. We denote this grid in graph notation as
$G_{\Omega} = (V_{\Omega},E_{\Omega})$. The entries of $A$ are then
given as 
\begin{equation*}
a_{i,j} = \left\{ \begin{array}{cl} \frac{1}{d_{out}(j)}, & \text{if\
  } (i,j) \in E_{\Omega}\\
  0, & \text{otherwise}, \end{array}\right.
\end{equation*} where $d_{out}(j)$ is the number of outgoing edges of
$j \in \Omega$. In Figure~\ref{fig:nr:scalar:uniformnet}, we illustrate
the two-dimensional uniform network problem.
In the tests we conduct, we use again full-coarsening, i.e., we choose
$\coarsevar$ to be every other grid point in both spatial dimensions. In order
to keep the overall invested work small, we consider to take only up
to two interpolatory points $\coarsevar_{i}$ per point $i \in
\finevar = \Omega \setminus \coarsevar$.
As the steady-state vector is known to be strictly positive, we
choose to use initially random, but positive test vectors for 
this first problem, and also in all following tests. More
precisely, we choose vectors with entries uniformly distributed in $[1,2]$.

In Table~\ref{tab:2DUniformNetwork}, we present results that use $6$
test vectors and a $V(2,2)$-cycle MLE step with $\omega$-Jacobi
smoother, $\omega = .7$. We report the number of iterations needed to
compute the steady-state vector $x$, such that
\begin{equation*}
  \norm{Bx} \leq 10^{-8}, \norm{x} = 1.
\end{equation*} In addition, the number of preconditioned GMRES 
iterations needed to achieve the same accuracy is reported, with the 
initial MLE setup cycle denoted by the subscript.
The coarsest grid in the experiments is always $5 \times 5$.
\begin{table}[!ht]
  \begin{center}
    \begin{tabular}{|c|c|c|c|c|cc|cc|}
      \hline
      $N$ & 17 & 33 & 65 & 129 \\
      \hline
      MLE & \AF{$10$} &  \AF{$9$} & \AF{$10$} & \AF{$11$}\\ 
      \hline
      pGMRES & \AF{$7_{V}$} & \AF{$8_{V}$} &
      \AF{$10_{V}$}  & \AF{$10_{V^{2}}$}  \\
      \hline
      pArnoldi & \AF{$7_{V}$} & \AF{$8_{V}$} &
      \AF{$10_{V}$}  & \AF{$9_{V^{2}}$}  \\
      \hline
    \end{tabular}
    \caption[BootAMG ($k=6, \eta=2$) MLE for uniform network problem]{Multilevel results for the two-dimensional uniform
      network model on an $N\times N$ grid. We report results to
      compute the steady-state vector $x$ to an accuracy of $10^{-8}$,
      using a $V(2,2)$-MLE cycle with $\omega$-Jacobi smoothing,
      $\omega=.7$. In addition we also report the number of iterations 
      pGMRES needs to achieve the same accuracy, where we denote the
      initial bootstrap setup in the subscript. 
      The sets $\mathcal{U}$ and $\mathcal{V}$ 
      consist of $6$ initially positive random vectors and coarsest grid
      eigenvectors, respectively.\label{tab:2DUniformNetwork}}
  \end{center}
\end{table} The operator complexity in these tests is bounded by
$1.6$. The test shows that the MLE scales with increasing problem-size, while
the number of preconditioned GMRES iterations grows slightly. However, one step of pGMRES is much cheaper than one
step of MLE.  Although we do not report results here, in general we may also consider restarting the MLE iterations when 
the solution is not found in a reasonable number of pGMRES iterations. 

\begin{figure}[!t]
  \centering 
  \scalebox{.75}{\input{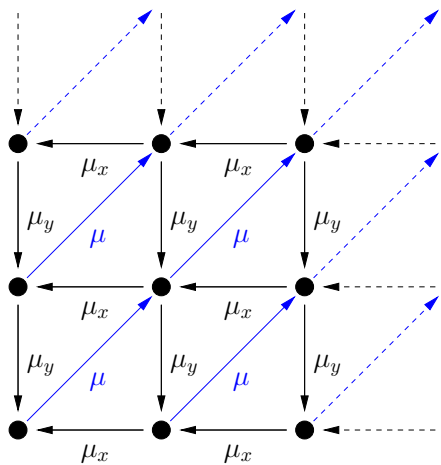}}
  \caption[Tandem-Queuing network model]{Tandem-Queuing Network with probability to advance $\mu$
    and queuing probabilities $\mu_{x}$ and $\mu_{y}$ on a
    two-dimensional equidistant grid.\label{fig:nr:scalar:tdnet}}
\end{figure}

\begin{figure}
  \begin{center}
      \includegraphics[width=\textwidth]{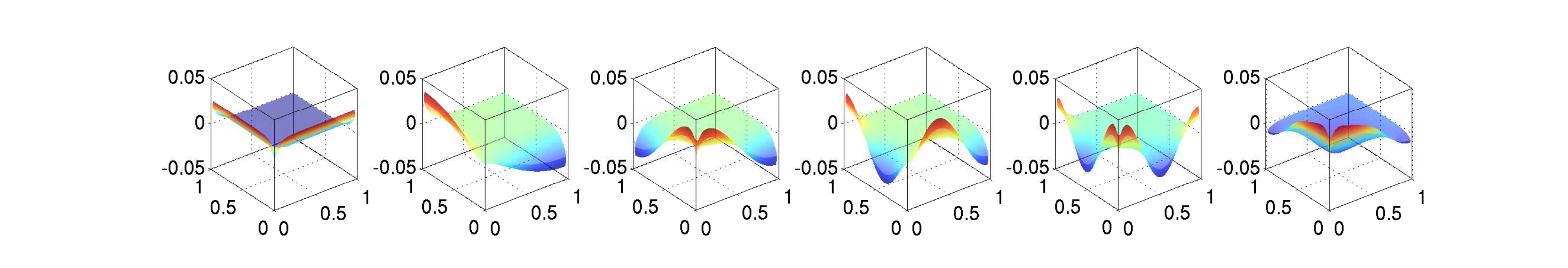}
  \end{center}
  \caption[Eigenvector approximations for the TQN problem]{Approximations to the eigenvectors corresponding to the
    smallest $6$ eigenvalues of the tandem queuing network problem on
    a $129 \times 129$ grid with $6$ levels upon convergence of the
    steady-state solution of the MLE method.\label{fig:tqnevc129}}
\end{figure}

The next Markov chain model we consider in our tests is a tandem
queuing network, illustrated in Figure~\ref{fig:nr:scalar:tdnet}. 
 The spectrum of this operator is complex, as we 
show in Figure~\ref{fig:specplotJ}. Again, we use
full-coarsening and a coarsest grid of $5\times 5$. We present 
our results in Table~\ref{tab:TandemQueuingNetwork}. In the
test we use the following probabilities,
\begin{equation*}
  \mu = \frac{11}{31}, \quad \mu_{x} = \frac{10}{31}, \quad \mu_{y} =
  \frac{10}{31}.
\end{equation*}
\begin{table}[!t]
  \begin{center}
    \begin{tabular}{|c|c|c|c|c|cc|cc|}
      \hline
      $N$ & 17 & 33 & 65 & 129 \\
      \hline
      MLE & \AF{$8$} & \AF{$8$} &
      \AF{${8}$} & \AF{${8}$} \\ 
      \hline
      pGMRES & \AF{${6_{V}}$} & \AF{$6_{V}$} & 
      \AF{$6_{V}$} & \AF{$7_{V}$} \\
      \hline
      pArnoldi & \AF{${6_{V}}$} & \AF{$6_{V}$} & 
      \AF{$6_{V}$} & \AF{$7_{V}$} \\
      \hline
    \end{tabular}
    \caption[BootAMG ($k=6, \eta=2$) MLE for TQN problem]{Multilevel results for the tandem queuing
      network model on an $N\times N$ grid. We report results to
      compute the steady-state vector $x$ to an accuracy of $10^{-8}$,
      using a $V(2,2)$-MLE cycle with $\omega$-Jacobi smoothing,
      $\omega=.7$. In addition we also report the number of iterations 
      pGMRES needs to achieve the same accuracy, where we denote the
      initial bootstrap setup in the subscript. 
      The sets $\mathcal{U}$ and $\mathcal{V}$ 
      consist of $6$ initially positive random vectors and coarsest grid
      eigenvectors, respectively.\label{tab:TandemQueuingNetwork}}
  \end{center}
\end{table} Again, we see that the MLE method converges rapidly to the
steady-state vector and also yields a very efficient preconditioner
for the GMRES method. We observe that the number of iterations does
not depend on the size of the problem. Note that the MLE method also
yields accurate approximations to the eigenvectors corresponding
to all $k$ smallest eigenvalues on the finest grid. 
In Figure~\ref{fig:tqnevc129}, we show the computed approximations 
and report in Table~\ref{tab:nr:accevctqn129} their
accuracy upon convergence of the steady-state vector.
\begin{table}[!ht]
  \begin{center}
    { \begin{tabular}{|c|c|c|c|c|c|c|c|c|}
      \hline
      $i$ & $1$ & $2$ & $3$ & $4$ & $5$ & $6$\\
      \hline       
      $\norm{v_i^{L} - v_i}$ & {\small$1.02E{-8}$} & {\small$9.04E{-3}$} &
      {\small$2.68E{-2}$} & {\small$2.84E{-2}$} &
      {\small$1.42E{-1}$} & {\small$3.77E{-1}$}\\
      \hline
    \end{tabular}}
    \caption[Accuracy of EV approximations for the TQN problem]{Accuracy of the eigenvectors corresponding to the
    smallest $6$ eigenvalues of the tandem queuing network problem on
    a $129 \times 129$ grid with $6$ levels upon convergence of the
    steady-state solution of the MLE method.}\label{tab:nr:accevctqn129}  
  \end{center}
\end{table} One should keep in mind that the results are intended to show the promise
of our approach 
rather than presenting an optimized
method in the bootstrap framework for this type of problems. We limit
our analysis to the statement that with minimal effort spent in adjusting
the parameters of the setup of the bootstrap approach we obtain 
scalable solvers. The optimization of the method by tuning the
parameters involved, e.g., relaxation parameter of the smoother, coarsening,
caliber, number of test vectors, weighting in the least squares
interpolation, number of relaxations in the setup and solution
process, is part of future research.

\begin{figure}[!ht]
  \subfigure[Random planar graph with $N=256$\label{fig:rpg256}]
{
\includegraphics[width=.45\textwidth]{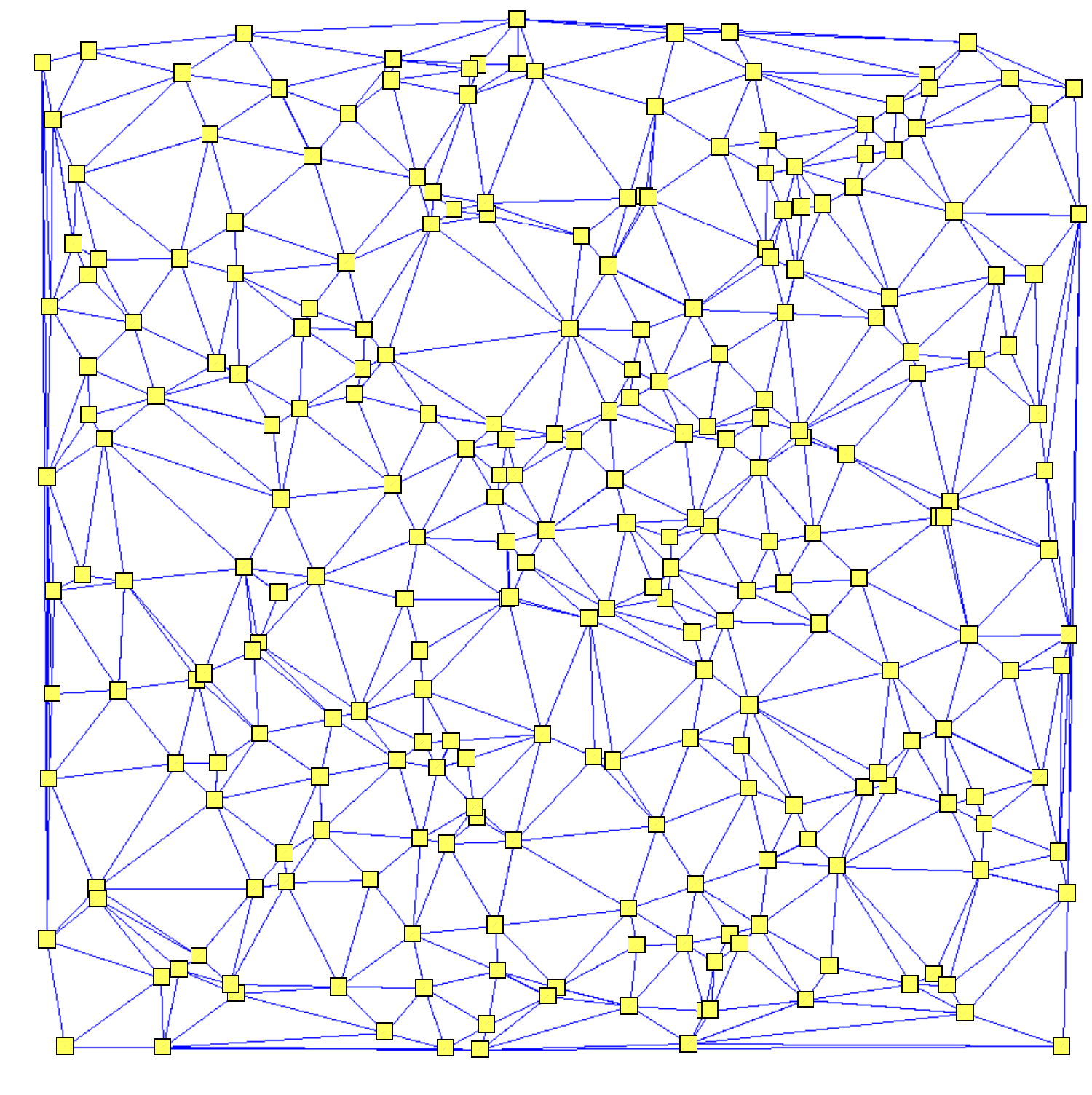}\hfill 
}
  \subfigure[Random planar graph with $N=1024$\label{fig:rpg1024}]{
\includegraphics[width=.45\textwidth]{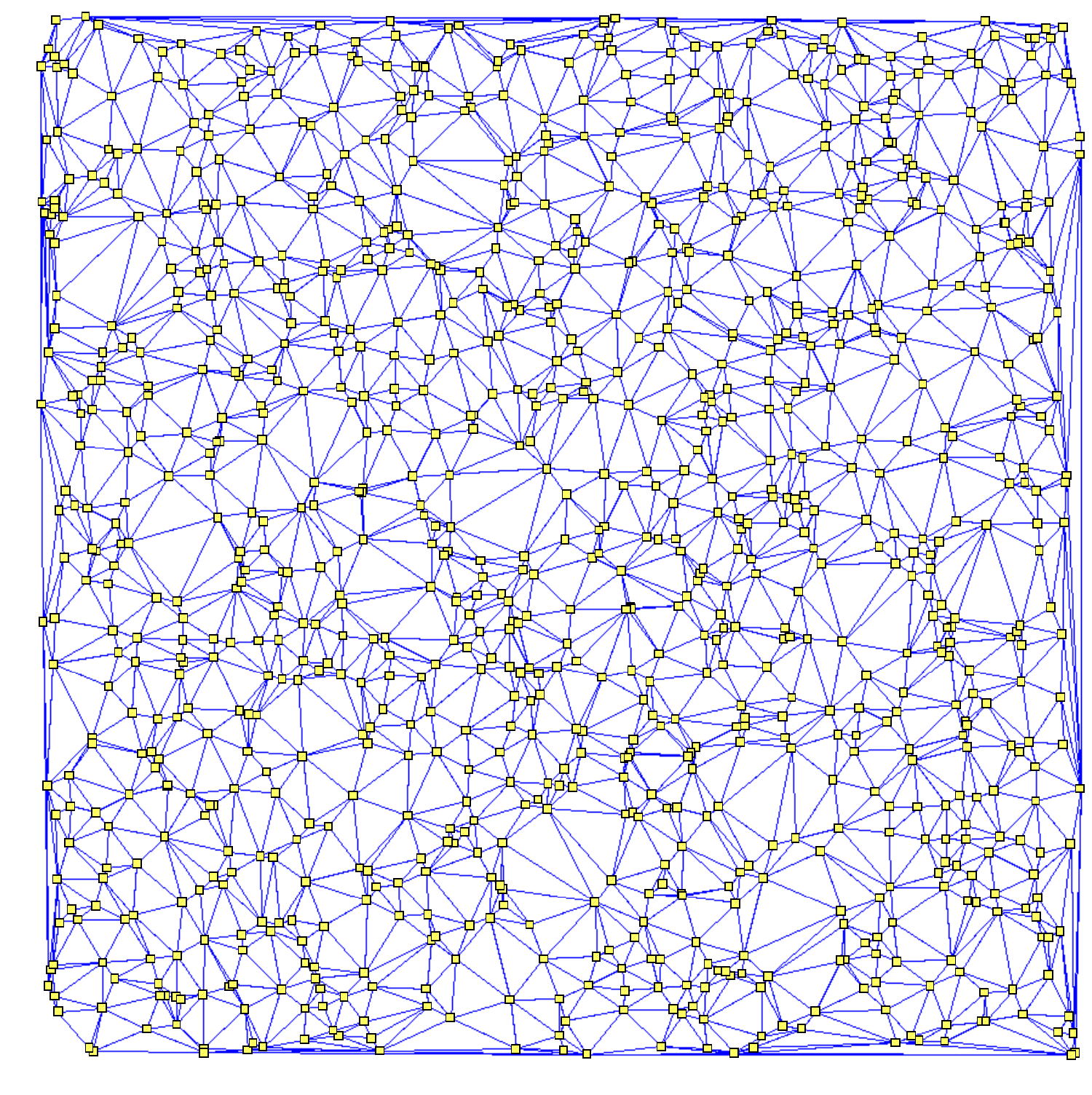}
}
  \caption[Unstructured planar graph network model]{Delaunay-triangulations of $N$ randomly chosen points in
  the unit square $[0,1]\times[0,1]$.\label{fig:randomplanargraph}}
\end{figure}

\begin{figure}[!ht]
  \subfigure[Random planar graph with $N=256$]
{
\includegraphics[width=.45\textwidth]{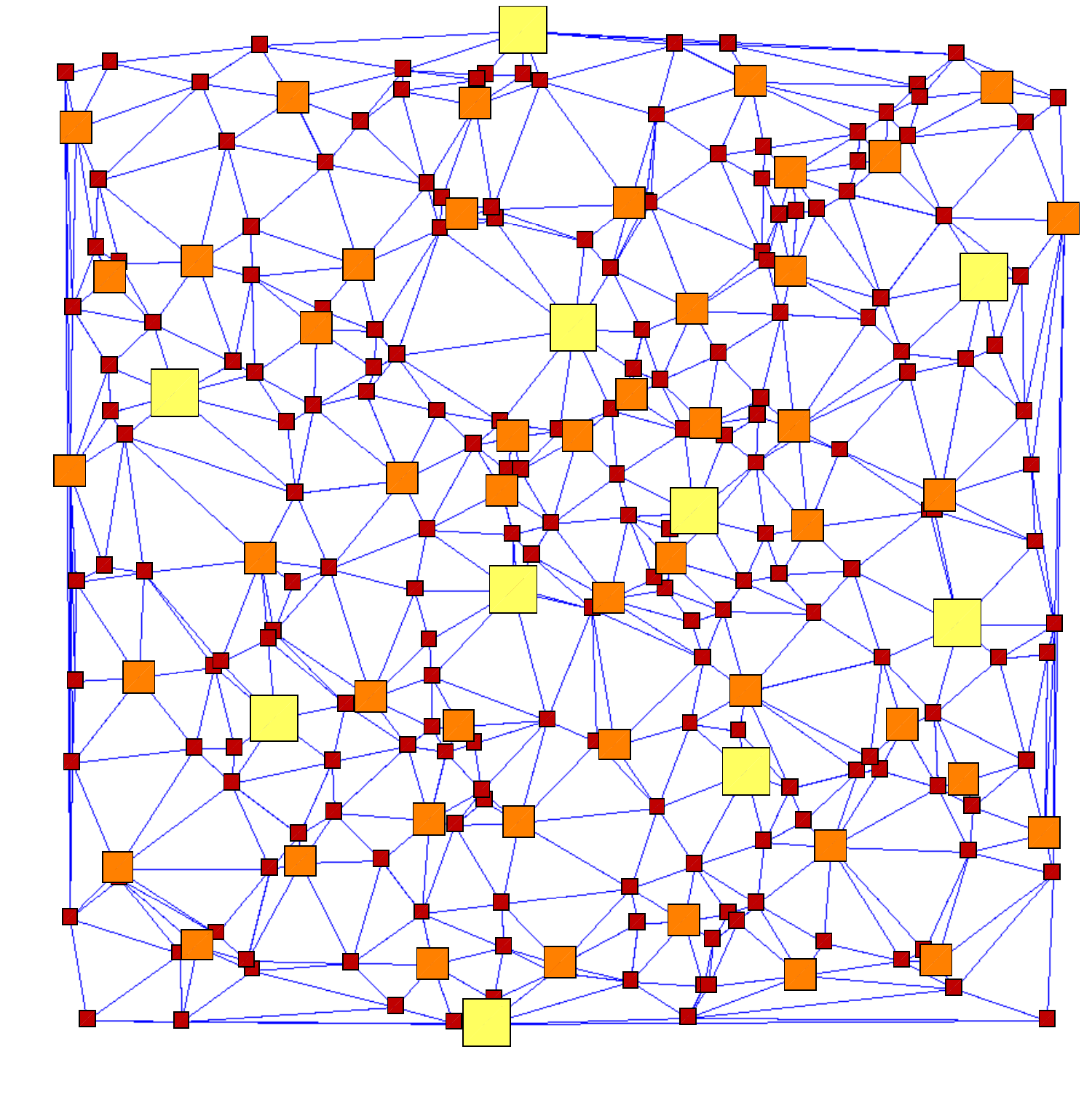}\hfill 
}
  \subfigure[Random planar graph with $N=1024$]
{
\includegraphics[width=.45\textwidth]{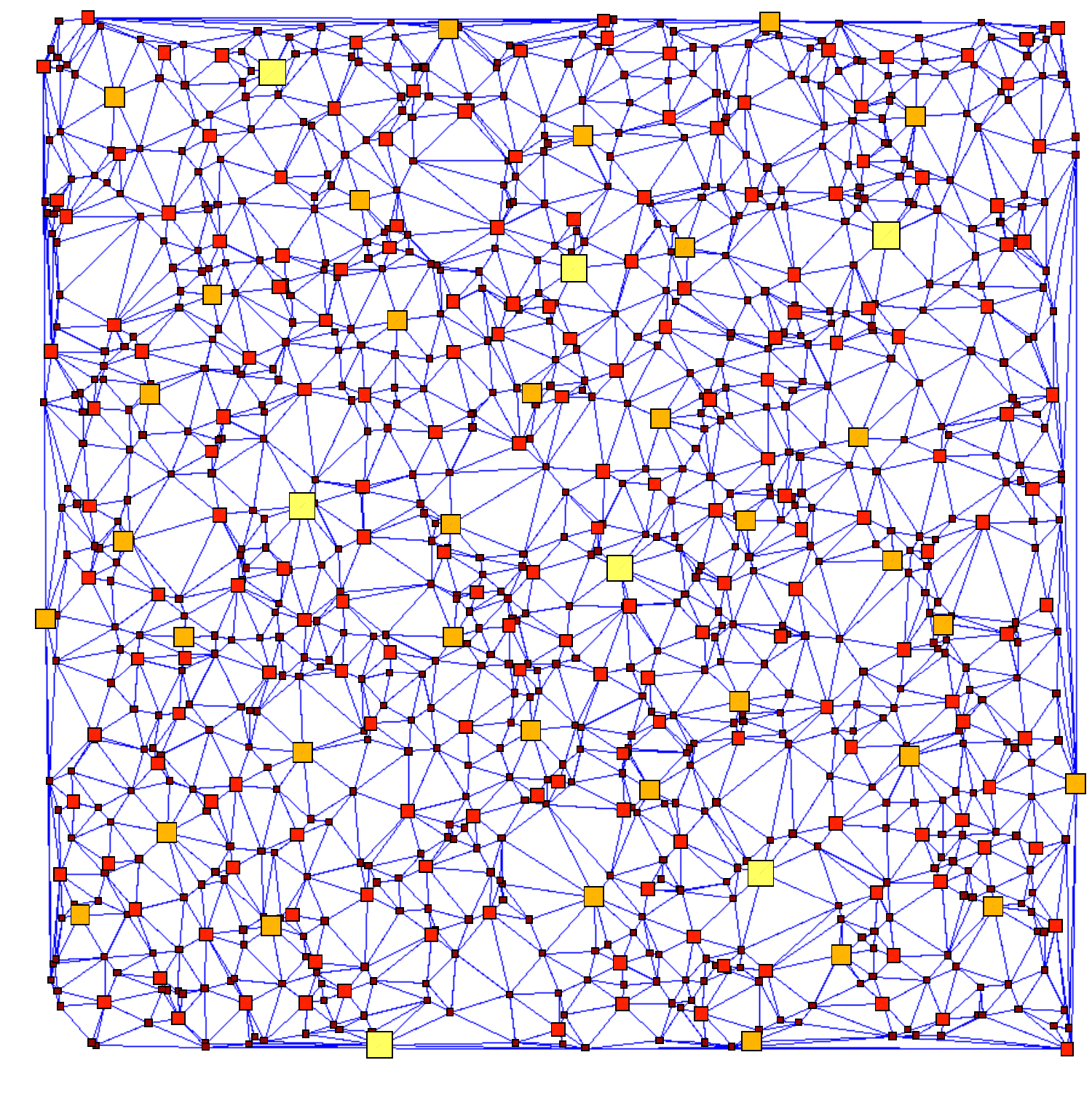} 
}
  \caption[Coarsening of UPG by compatible relaxation]{Coarsening of Delaunay-triangulations of $N$ points
  randomly scattered  in the unit square
  $[0,1]\times[0,1]$ shown in   Figure~\ref{fig:randomplanargraph}
  using compatible relaxation.\label{fig:randomplanargraphcoarse}} 
\end{figure}

The last test we consider corresponds to  a triangulation of $N$
randomly chosen points in $[0,1]\times[0,1]$. The transition
probabilities in the network are then given by \AF{the inverse of} the number of outgoing
edges at each point, similar to the uniform network. 
In Figure~\ref{fig:randomplanargraph}, 
we show two examples of such networks, one with $N=256$ and one with
$N=1024$ points. 
Due to the fact that the corresponding graphs of this model are
planar, we call this model \emph{unstructured planar graph} model.
As there is no natural way to define the set of coarse variables
$\coarsevar$ for this problem we use compatible relaxation,
introduced in section~\ref{sec:CompRel} to define the splitting of variables
$\Omega = \finevar \cup \coarsevar$. 
In Figure~\ref{fig:randomplanargraphcoarse}, a resulting
coarsening is presented. The size of each individual point represents
how many grids it appears on. 
In Table~\ref{tab:unstructuredplanarG}, we report results of our
overall algorithm for unstructured planar graph models for a
set of varying graphs. 
\begin{table}[!ht]
  \begin{center}
    \begin{tabular}{|c|c|c|c|c|cc|cc|}
      \hline
      $N$ & 256 & 512 & 1024 & 2048 \\
      \hline
      MLE & \AF{$15$} & \AF{$20$} & \AF{$20$} & \AF{$20$}  \\ 
      \hline
      pGMRES & \AF{$8_{V}$}  & \AF{$10_{V}$}  &
      \AF{$10_{V}$} &  \AF{$11_{V}$} \\
      \hline
      pArnoldi & \AF{$8_{V}$}  & \AF{$10_{V}$}  &
      \AF{$10_{V}$} &  \AF{$11_{V}$} \\
      \hline
    \end{tabular}
    \caption[BootAMG ($k=6, \eta=2$) MLE for UPG problem]{Multilevel
      \AF{($3$ level)} results for the unstructured planar graph
      model on a grid with $N$ points randomly scattered in the unit square. We report results to 
      compute the steady-state vector $x$ to an accuracy of $10^{-8}$,
      using a $V(2,2)$-MLE cycle with $\omega$-Jacobi smoothing,
      $\omega=.7$. In addition we also report the number of iterations 
      pGMRES needs to achieve the same accuracy, where we denote the
      initial bootstrap setup in the subscript. 
      The sets $\mathcal{U}$ and $\mathcal{V}$ 
      consist of $6$ initially positive random vectors and coarsest grid
      eigenvectors, respectively. \label{tab:unstructuredplanarG}
     }
  \end{center}
\end{table} Even for this unstructured graph network we obtain a fast
converging method with our MLE approach.  
The mild variation in the results when increasing the problem-size might be
caused by the fact that by increasing the problem-size 
the nature of the problem
changes in the sense that the average number of 
outgoing edges of grid points
increases. That is, it is not necessarily clear whether two unstructured
graphs of different sizes are comparable.

\section{Conclusions}\label{sec:Conclusions}
The proposed BootAMG Multilevel setup algorithm produces an effective multilevel eigensolver for the 
Markov-Chain test problems we considered.  In general, we have proposed and implemented two important new ideas:
First, the use of a coarse-level eigensolve and the resulting
multilevel hierarchy to improve a given approximation of the state
vector. Second the use of BootAMG preconditioned GMRES
steps to further accelerate this eigensolver.  
Both ideas, separately or combined can be incorporated
into  any given multilevel method used for  solving Markov Chain problems or other problems targeting smooth eigenvectors.
 The accurate representation of the near-kernel of the fine-level system on coarser levels, that results from our BootAMG setup, yields a very effective preconditioner to GMRES \AF{as well as for the Arnoldi method}.   
An additional benefit of our proposed method over other existing multilevel method for Markov Chains is that we do not require
any special processing of the coarse-level systems to ensure that
stochastic properties of the fine-level system are maintained there. 
We mention that the developed approach is not restricted to Markov Chain problems; it
can be applied to other eigenvalue problems.  

\bibliographystyle{siam}


\newcommand{\etalchar}[1]{$^{#1}$}
\def\cdprime{$''$}

\end{document}